\documentclass[11pt]{article}
\usepackage{amssymb,amsfonts,amsmath}

\begin{document}
\renewcommand{\hat}{\widehat}

\newcommand\hot{\hat\otimes}
\newcommand\DDD{\mathfrak D}
\newcommand\id{\operatorname{id}}

\def\wpx#1{{:}\,x^{\ot #1}\,{:}}
\newcommand\indlim{\operatornamewithlimits{ind\,lim}}
\newcommand\R{\mathbb R}
\newcommand\N{\mathbb N}
\newcommand\Z{\mathbb Z}
\newcommand\Zp{\mathbb Z_+}
\newcommand\C{{\mathbb C}}
\newcommand\SP{(S_{\text{\rom P}})}
\newcommand\EP{(\mathcal E_{\text{\rom P}})}

\newcommand\muP{\mu_{\text{\rom{P}}}}
\newcommand\LP{(L^2_{\text{\rom P}})}
\newcommand\rom[1]{\textup{#1}}
\newcommand{\vs}{\vspace{0.3cm}}
\newcommand{\lth}{$\lambda,\theta$}
\newcommand{\lt}{_{\lambda,\theta}}
\newcommand{\di}{\partial}
\newcommand{\dis}{\partial^*}
\newcommand{\be}[1]{\begin{equation}\label{#1}}
\newcommand{\bea}{\begin{eqnarray}}
\newcommand{\eea}{\end{eqnarray}}
\newcommand{\beas}{\begin{eqnarray*}}
\newcommand{\eeas}{\end{eqnarray*}}
\newcommand{\n}{{\mathbb N}}
\newcommand{\s}{{\mathcal S}}
\newcommand{\ds}{&\displaystyle}
\newcommand{\cd}{\cdot}
\newcommand{\e}{{\mathcal E}}
\newcommand{\prlim}{\mathop{\mathrm{proj\, lim}}\limits_{p\to\infty}}
\newcommand{\ho}{\hat \otimes}
\newcommand{\ot}{\otimes}
\newcommand{\x}[1]{:\!x^{\ot #1}\!:}
\newcommand{\y}[1]{:\!y^{\ot #1}\!:}
\newcommand{\yj}[1]{:\!y_j^{\ot #1}\!:}
\newcommand{\A}{{\mathcal A}}
\newcommand{\dm}{\diamondsuit}
\newcommand{\p}{{\mathcal P}}
\newcommand{\bt}{\bold t}
\newcommand{\bs}{\bold s}
\newcommand{\bu}{\bold u}
\newcommand{\wexp}[1]{:\!e^{<x,#1>}\!:}

\newtheorem{theorem}{Theorem}[section]
\newtheorem{corollary}{Corollary}[section]
\newtheorem{prop}{Proposition}[section]
\newtheorem{lemma}{Lemma}[section]

\makeatletter\@addtoreset{equation}{section}\makeatother
\def\theequation{\arabic{section}.\arabic{equation}}

\title{\Large  A NOTE ON SPACES OF TEST AND GENERALIZED FUNCTIONS\\ \Large
OF POISSON WHITE NOISE}

\author{Eugene W. Lytvynov\thanks{Supported in part by  AMS grant}}
\date{}
\maketitle
\begin{abstract} The paper is devoted to construction and investigation
of some riggings of the $L^2$-space of Poisson white noise. A particular attention
is paid to the existence of a continuous version of a function from a
test space,
and to the property of an algebraic structure under pointwise multiplication
of functions from a test space.\end{abstract}


\section{Introduction}
In the works \cite{It,ItKu}, started was a study of test and generalized
functions defined on the Schwartz space of tempered distributions $\s'(\R)$
the dual pairing of which is determined by the inner product of the $L^2$-space
$\LP=L^2(S'(\R),d\muP)$, where $\muP$ is the measure of Poisson white noise.
Following the construction of the space of Hida distributions in Gaussian analysis
(e.g., \cite{KuTaI,BeKo,HKPS,Ob2}),
Ito and Kubo \cite{ItKu} introduced the triple $\SP^*\supset\LP\supset\SP$.
However, the two following important problems remained open: 1) Does the space $
\SP$ consist of continuous functions, or, which is ``almost'' equivalent,
do the delta functions belong to $\SP^*$? 2) Is the space $\SP$ an algebra
under pointwise multiplication of functions?

\renewcommand{\kappa}{\varkappa}

In this note, we will construct a whole scale of test spaces $\SP^\kappa$,
$\kappa\ge0$, such that $\SP^{\kappa_1}\subset\SP^{\kappa_2}$ if $\kappa_1>
\kappa_2$, and of their dual spaces $\SP^{-\kappa}$ with respect to $\LP$.
For $\kappa=0$, $\SP^0=\SP$, so that $\SP^{-0}=\SP^*$. The idea of construction
of these spaces comes  from the corresponding constructions in
Gaussian analysis \cite{KS,KLS}.

The main results of the paper are the following: 1) The space $\SP^{1}$
consists of continuous functions, and for each $\SP^\kappa$ with $\kappa<1$
this is not the case. 2) The space $\SP^1$ (and even each space $\SP^\kappa$
with $\kappa>1$) is an algebra under pointwise multiplication, and,
moreover, an estimate of Hilbert norms of a product of two functions
from $\SP^1$ is given, this estimate being rather analogous to the
estimates which hold
for the Hida test space in Gaussian analysis, e.g.,
\cite{KuTaII,Yokoi1,Yokoi2,HKPS}.

It should be noted that there is another approach to developing analysis on
non-Gaussian spaces, based on the use of a system of Appell polynomials and
its dual (biorthogonal) system \cite{ADKS,KSWY}.  By using Theorem~3.2 of the present paper,
one can prove that the application of the  biorthogonal analysis developed
in \cite{KSWY} to the Poisson measure leads, in fact, to the same triple
\be{1.1}\SP^{-1}\supset \LP\supset \SP^1,\end{equation}
see \cite{Us,LyUs}.
So, it seems that the triple \eqref{1.1} will play the role analogous to the
Hida triple in Gaussian analysis, and one will use either the orthogonal or
biorthogonal  methods, or combine them, depending on a specific problem under
study.

\section{Setup for Poisson white noise calculus}
In this section, we will construct the above mentioned scale of test and dual
spaces. To this end, we present below some results of the works
\cite{It,ItKu,Ly}, see also \cite{Su}.

Let $T$ be a separable, topological space and $\nu$ a $\sigma$-finite,
 non-atomic measure
defined on the Borel $\sigma$-algebra ${\mathcal B}(T)$. We consider
a standard (Gelfand)  triple of spaces \cite{GeVi}
\be{2.1}
\e'=\indlim_{p\to\infty}E_{-p}\supset L^2(T,\nu)=E_0\supset\projlim_{p\to
\infty}E_p=\e,\end{equation}
where
$\{E_p\, |\, p\ge0\}$ is a sequence of real, compatible, separable, Hilbert
spaces such that, for any $p>q\ge0$, $E_p$ is topologically---that is,
densely and continuously---embedded into $E_q$ and
\be{2.2} |\xi|_p \ge |\xi|_q,\qquad \xi\in E_p,\end{equation}
where $|\cdot|_p$ denotes the $E_p$ norm; $E_{-p}$ is the dual of $E_p$
with respect to zero space $E_0$, so that $\e'$ is the dual of $\e$.

One  makes the following assumptions about the space
$\e$:\vs

(A.1) Every element of $E_1$ has a version continuous on $T$, and for
each $t\in T$, $\delta_t$---the delta function at $t$---belongs to $E_{-1}$.
Moreover, $\delta\colon t\to\delta_t$ is a continuous mapping of $T$ into $E_{-1}$
and
$$\|\delta\|^2\equiv \int_T|\delta_t|_{-1}^2 \, d\nu(t)<\infty,$$
so that \cite{KuTaI} the embedding operator $E_1\hookrightarrow E_0$ is of
Hilbert--Schmidt type. Its Hilbert--Schmidt norm is less than 1.

(A.2) The mapping $\delta$ satisfies
$$\|\delta\|_\infty \equiv \int_T|\delta_t|_{-1}\, d\nu(t)+ \sup_{t\in T}
|\delta_t|_{-1}<\infty,$$
so that $E_1\subset L^1(T,\nu)\cap L^\infty(T,\nu)$.

(A.3) For any $\xi\in\e$ and $p\ge0$,
$$\rho|\xi|_{p+1}\ge|\xi|_p$$
with a fixed $\rho\in(0,1)$.

(A.4) The ``diagonalization'' operator ${\mathfrak D}:\e^{\ot 2}\to\e$,
$\displaystyle\e^{\ot2}=\projlim_{p\to\infty} E_p^{\ot2}$, given by
\be{star}({\mathfrak D}f^{(2)})(t)=f^{(2)}(t,t),\qquad f^{(2)}\in\e^{\ot2},\
t\in T,
\end{equation}
is continuous, and moreover, for any $p\ge1$,
$$|{\mathfrak D}f^{(2)}|_p\le C_p|f^{(2)}\|_p,\qquad C_p>0,$$
$|\cdot|_p$ standing also for the norm of each space $E_p^{\otimes n}$,
so that $\e$ is an algebra under pointwise multiplication of functions and
$$|\xi\eta|_p\le C_p|\xi|_p|\eta|_p,\qquad p\ge1.  $$

(A.5) The set of the functions $\xi\in\e$ whose support is  of
finite $\nu$ measure is dense in $\e$.\vs

Let us construct an example of such a triple. Fix the sequence $(e_j)_{j=0}
^\infty$ of the Hermite functions on $\R$:
$$ e_j=e_j(t)=(\sqrt\pi\,2^jj!)^{-1/2}(-1)^je^{t^2/2}(d/dt)^je^{-t^2}.$$
For each $p\ge1$, define $\s_p(\R)$ to be the real Hilbert space spanned
by the orthonormal basis $\big(e_j(2j+2)^{-p}\big)_{j=0}^\infty$, and let
$\s_p(\R^d)\equiv \s_p(\R)^{\ot d}$.
Considered as a subset of $L^2(\R^d)$, every space $\s_p(\R^d)$ coincides with
the domain of the operator $(H^{\ot d})^p$, where $H^{\ot d}$ is the harmonic
oscillator in $L^2(\R^d)$:
$H^{\ot d}=-\sum_{i=1}^d\big(\frac{d}{dt_i}\big)^2+
\sum_{i=1}^d t_i^2+1$.
As  well known, $\s(\R^d)=\projlim_{p\to\infty}
\s_p(\R^d)$ is the Schwartz space of rapidly decreasing functions
on $\R^d$. Denote by $\s_{-p}(\R^d)$ the dual of $\s_p(\R^d)$. Then,
$\s_1(\R^d)$ consists of continuous functions and
$\R^d\ni t\to \delta_t\in\s_{-1}(\R^d)$
is a continuous mapping such that
$\sup_{t\in\R^d}\|\delta_t\|_{\s_{-1}(\R^d)}<\infty$.
The assumptions (A.3), (A.4), and (A.5) are satisfied for $\s_p(\R^d)$'s.
Let now $\nu$ be a $\sigma$-finite, non-atomic, Borel, regular,
measure on $\R^d$. Suppose also
that, for some $\varepsilon\ge0$,
\be{2.3}\int_{\R^d}\big( \|\delta_t\|_{\s_{-1-\varepsilon}(\R^d)}^2+\|\delta_t\|_{\s_{-1-
\varepsilon}(\R^d)}\big)\, d\nu(t)<\infty\end{equation}
(for the Lebesgue measure, this holds when $\varepsilon=0$). Making use
of the evident estimate
\beas
\ds \|\xi\|^2_{L^2(\R^d,\nu)}=\int_{\R^d}|\xi(t)|^2\,d\nu(t)\le\|\xi\|_{\s_{1+
\varepsilon}(\R^d)}\int_{\R^d}\|\delta_t\|_{\s_{-1-\varepsilon}(\R^d)}^2\,d\nu(t)
,&\\
\ds \xi\in\s_{1+\varepsilon}(\R^d),&\end{eqnarray*}
we conclude that $\s_{1+\varepsilon}(\R^d)$ is continuously embedded into
$L^2(\R^d,\nu)$. Moreover, by (\ref{2.3}) and \cite{KuTaI}, the embedding operator
$O_{1+\varepsilon}:\s_{1+\varepsilon}(\R^d)\hookrightarrow L^2(\R^d,\nu)$
is of Hilbert--Schmidt type (note that, by passing to an equivalent system of
norms, one can always make its Hilbert--Schmidt norm  less than 1). Then
because of the regularity and $\sigma$-finiteness of $\nu$, $\s(\R^d)$
is a dense subset of $L^2(\R^d,\nu)$. At last, for each $p\ge1$, define
$E_p$ to be the Hilbert factor space $\s_{p+\varepsilon}/\ker O_{p+
\varepsilon}$, where
$O_{p+\varepsilon}\colon\s_{p+\varepsilon}(\R^d)\hookrightarrow L^2(\R^d,\nu)$
is an embedding operator. By \cite{BeKo}, Ch.~5, Sect.~5, subsec.~1,
$\{L^2(T,\nu), \, E_p\, |\, p\ge1\}$ is a sequence of compatible, Hilbert
spaces,  where $T$ denotes  the support of $\nu$. Thus, we get the desired triple
\be{2.4}
\s'(T)=\e'=\indlim_{p\to\infty} E_{-p}\supset L^2(T,\nu)\supset \prlim E_p=\e=\s(T).
\end{equation}
Note that  the spaces $E_p$, $p\ge1$, are completely determined
by the set $T$, and $\e$ is actually the test Schwartz space on $T$,
$\s(T)$, respectively $\e'$ is the  Schwartz space of tempered distributions
on $T$, which is the dual of $\s(T)$  with respect to zero space $L^2(T,\nu)$.
Note also that, in case of a bounded $T$, $\s(T)={\mathcal D}(T)$ is the space
of infinitely differentiable functions on $T$.

Given a real Hilbert space ${\mathcal H}$ and $\kappa\in\R$, a weighted Fock space
$\Gamma_\kappa({\mathcal H})$ is defined by
$$\Gamma_\kappa({\mathcal H})\equiv  \bigoplus_{n=0}^\infty {\mathcal H}_{\C}
^{\ho n}(n!)^{1+\kappa},$$
where the symbol $\hot$ denotes the symmetric tensor product, the index
$\C$ stands for  complexification of a real space, $\mathcal H_\C^{\hot 0}=\C$,
$0!=1$.
Particularly, $\Gamma_0({\mathcal
H})=\Gamma({\mathcal H})$ is the usual Fock space over ${\mathcal H}$.

By using (\ref{2.1}) and (\ref{2.2}), we construct, for each  $\kappa\ge0$, the following
standard triple \cite{KuTaI,BeKo,HKPS,KS,KLS}
\be{2.5}
\Gamma_{-\kappa}(\e')=\indlim_{p\to\infty} \Gamma_{-\kappa}(E_{-p})\supset\Gamma(E_0)\supset
\prlim \Gamma_\kappa(E_p)=\Gamma_\kappa(\e).\end{equation}
By (A.1), the embedding operator
$\Gamma_\kappa(E_1)\hookrightarrow \Gamma(E_0)$ is of Hilbert--Schmidt type,
see, e.g., \cite{BeKo,HKPS}.

We will use also the following triple
\be{2.6}
\Gamma_{\mbox{\scriptsize fin}}(\e)^*\supset \Gamma(E_0)
\supset\Gamma_{\mbox{\scriptsize fin}}(\e).\end{equation}
Here, $\Gamma_{\mbox{\scriptsize fin}}(\e)$ is the topological direct
sum  of the spaces
$$ \e_{\C}^{\ho n}\equiv\prlim E_{p,\C}^{\ho n},\qquad
n\in\Z_+=\{0,1,2,\dotsc\},$$
so that $\Gamma_{\mbox{\scriptsize fin}}(\e)$ consists of finite sequences
$(f^{(n)})_{n=0}^\infty$, $f^{(n)}\in\e_{\C}^{\ho n}$. The convergence
in $\Gamma_{\mbox{\scriptsize fin}}(\e)$ is equivalent to the uniform finiteness
and the coordinate-wise convergence in $\e_{\C}^{\ho n}$. The
$\Gamma_{\mbox{\scriptsize fin}}(\e)^*$ is the dual of
$\Gamma_{\mbox{\scriptsize fin}}(\e)$ with respect to $\Gamma(E_0)$.
It consists of all the sequences of the form $(F^{(n)})_{n=0}^\infty$,
$F^{(n)}\in\e_{\C}^{\prime\, \ho n}$, where
$$ \e_{\C}^{\prime\, \ho n}\equiv
\indlim_{p\to\infty} E_{-p,\C}^{\ho n},
\qquad n\in\Z_+.$$
The convergence in $\Gamma_{\mbox{\scriptsize fin}}(\e)^*$ is the
coordinate-wise convergence in $\e_{\C}^{\prime\, \ho n}$.

Now, endow $\e'$ with the strong dual topology, and define on the Borel
$\sigma$-algebra $\mathcal B(\e')$ the probability measure $\muP$ by its
Fourier transform
$$\int_{\e'}e^{i<x,\xi>}\,d\muP(x)=\exp\bigg[\int_T\big(
e^{i\xi(t)}-1\big)\,d\nu(t)\bigg],\qquad \xi\in \e,$$
$\muP$ is called the measure of Poisson white noise on $T$ with intensity $\nu$.
Here, $<\cdot\,,\cdot>$ stands for dualization between the space $\e_\C^{\prime\,
\hot n}$ and $\e_\C^{\hot n}$ for each $n$, which is supposed to be linear in both
dots.

For any $x\in\e'$, the Poisson Wick power $\wpx{n}\in\e^{\hot n}$,
$n\in\Z_+$, is defined by
the recursion relation \be{recur}\begin{gathered}
\wpx{0}=1,\qquad\wpx{1}=x-1,\\
<\wpx{(n+1)},f^{(n+1)}>=<\wpx{n}\hot\wpx{1},f^{(n+1)}>\\
-n<\wpx{n},\mathfrak D^{(n+1)}f^{(n+1)}>-n<\wpx{(n-1)}\hot\tau,f^{(n+1)}>,\\
f^{(n+1)}\in \e^{\hot (n+1)},\ n\in\N,\end{gathered}\end{equation}
where $\mathfrak D^{(n+1)}\colon \e^{\hot(n+1)}\to\e^{\hot n}$ is the
continuous operator  given by
\begin{gather} \DDD^{(2)}=\DDD,\notag\\
\DDD^{(n+1)}=\id^{\ot(n-1)}\ot\,\DDD+\id^{\ot(n-2)}\ot\,\DDD\ot\id+\dots+\DDD\ot\id^{
\ot(n-1)},\qquad n\ge2,\label{starr}\end{gather}
 where id is the identity operator, $\DDD$ is defined by \eqref{star},
 and $\tau$ is
an element of $\e^{\hot2}$ such that
$$<\tau,f^{(2)}>=\int_T\big(\DDD f^{(2)}\big)(t)\, d\nu(t),\qquad
f^{(2)}\in \e^{\hot2}.$$

\noindent{\it Remark\/} 2.1. In order to distinguish between the Gaussian and Poisson Wick
powers, it would be better to denote the latter by $\wpx{n}_{\text{P}}$. But
since  the Gaussian Wick powers do not appear in this note, we do not
use such a notation.\vs

Evidently, for any $f^{(n)}\in \e_\C^{\hot n}$,
the dualization $<\wpx{n},f^{(n)}>$ is well defined and is a continuous function of
$x\in\e'$, which is called the Wick monomial with kernel $f^{(n)}$. Then, for
each $f^{(n)}\in\hat L^2(T^n,\nu^n)=\big(L^2(T,\nu)\big)_\C^{\hot n}$, we define
a function $<\wpx{n},f^{(n)}>$ as an element of the space $\LP=L^2(\e',d\muP)$
that is the $\LP$-limit of an arbitrary sequence $\big(<\wpx{n},f^{(n)}_j>
\big)_{j=0}^\infty$ such that $f^{(n)}_j\in\e_\C^{\hot n}$ and $f^{(n)}_j\to
f^{(n)}$ as $j\to\infty$ in $\hat L^2(T^n,\nu^n)$.

Next, for any $f\in L^2(T,\nu)\cap L^1(T,\nu)$, we  put
$$<x,f>\equiv <\wpx 1,f>+\int_Tf(t)\,d\nu(t)\in\LP.$$
Hence, for any set $\alpha\subset T$ of finite $\nu$ measure,
we can put $X_\alpha=X_\alpha(x)=<x,\chi_\alpha>\in\LP$, where $\chi_\alpha$
is the indicator of $\alpha$.
Then, $X_\alpha$ is the Poisson random measure on $T$ with intensity $\nu$,
 i.e., for any $n\in\N$ and for arbitrary disjoint sets $\alpha_1,\dots,\alpha_n
 \in\mathcal B(T)$, the random variables $X_{\alpha_1},\dots,X_{\alpha_n}$ are
 independent, and for each $\alpha$ $X_\alpha$ has the Poisson distribution with mean
 $\nu(\alpha)$. Thus, $\tilde X_\alpha=X_\alpha-\nu(\alpha)=<\wpx{1},\chi_\alpha>$
 is the centered Poisson random measure.

 \begin{prop}For each $f^{(n)}\in\hat L^2(T^n,\nu^n)$\rom,
 \be{2.7}<\wpx{n},f^{(n)}>=\int_{T^n}f^{(n)}(t_1,\dots,t_n)\,d\tilde X_{t_1}
 \dotsm d\tilde X_{t_n},\end{equation}
 where the right hand side of \eqref{2.7} is the $n$-fold Wiener--It\^o
 integral of $f^{(n)}$ by the centered Poisson random measure
 $\tilde X_\alpha$\rom.\end{prop}

Since a  centered Poisson random measure has the chaotic representation
property, Proposition~2.1 yields

\begin{theorem}
The following mapping is a unitary\rom:
\[\Gamma(L^2(T,\nu))\ni f=(f^{(n)})_{n=0}^\infty\to If=(If)(x)=\sum_{n=0}^
\infty
<\wpx{n},f^{(n)}>\in\LP.\]
\end{theorem}

So, we are able now to construct different riggings of $\LP$. We have only to apply the
unitary $I$ (or its extension by continuity) to the riggings of the Fock space
$\Gamma(L^2(T,\nu))$ and get corresponding riggings of $\LP$.

First, we note that
\begin{prop} We have
$$ I\big(\Gamma_{\mbox{\rom{\scriptsize fin}}}(\e)\big)={\mathcal P}(\e'),$$
where ${\mathcal P}(\e')$ is the set of all continuous polynomials on $\e'$---%
that is, the set of all complex functions on $\e'$ of the form
$\sum_{i=0}^n <x^{\ot n},f^{(n)}>$\rom, $f^{(n)}\in\e_{\C}^{\ho n}$\rom,
 $n\in\Z_+.$
\end{prop}

Thus, the application of $I$ to the rigging \eqref{2.6} gives
$$\mathcal P(\e')^*\supset\LP\supset\mathcal P(\e').$$

\noindent {\it Remark\/} 2.2. The topology of the nuclear space $\mathcal P(\e')$ is
supposed to be that induced from $\Gamma_{\text{fin}}(\e)$ by the isomorphism $I$.
One can also define a nuclear topology on $\mathcal P(\e')$ from that of $\Gamma_{
\text{fin}}(\e)$ by using the following natural isomorphism \cite{BeKo,KSWY}
\[\Gamma_{\text{fin}}\ni f=(f^{(n)})_{n=0}^\infty\to \mathcal U f=(\mathcal Uf)(x)=
\sum_{n=0}^\infty<x^{\ot n},f^{(n)}>. \]
But, in fact, these two topologies coincide.

Thus, every generalized function $\Phi$ from the biggest (in a sense) space
$\mathcal P(\e')^*$ can be represented in the form
\[\Phi=\Phi(x)=\sum_{n=0}^\infty<\wpx{n},F^{(n)}> ,\qquad F^{(n)}\in\e_\C^{\prime\,\hot
n},\]
and the dual pairing between $\Phi$  and a continuous polynomial
$$\phi(x)=\sum_{n=0}^\infty<\wpx{n},f^{(n)}>,\qquad f^{(n)}\in\e_
\C^{\hot n},$$
is given by
\[\ll\Phi,\phi\gg=\sum_{n=0}^\infty<\overline{F^{(n)}},f^{(n)}>n!,\]
$\overline{F^{(n)}}$ denoting the complex conjugate of $F^{(n)}$.

Next, by  applying  $I$  to (\ref{2.5}), we get
\be{2.88}
\EP^{-\kappa}=\indlim_{p\to\infty}\EP^{-\kappa}_{-p}\supset \LP\supset\prlim
\EP_p^\kappa=\EP^\kappa,\qquad \kappa\ge0\end{equation}
(we are using  natural notations for all the images of the spaces from
(\ref{2.5})). The norm of any Hilbert space
$\EP_{\sharp p}^{\sharp \kappa}$, $p\ge0$, $\kappa\ge0$, $\sharp\in\{
+,-\}$,  will be denoted by $\|\cdot\|_{\sharp\kappa,\sharp p}$.
The triple \eqref{2.88}  with $\kappa=0$ was investigated in \cite{ItKu}.

Let us consider two important examples of distributions.
The first one is a Poisson white noise monomial:
\begin{equation}\label{2.8}:\!X_{t_1}'\dotsm X_{t_n}'\!:=:\! x(t_1)\dotsm x(t_n)\!:
\equiv
<\wpx{n},\delta_{t_1}\ho\dotsb\ho\delta_{t_n}>,\qquad
t_1,\dots,t_n\in T.\end{equation}

By (A.1), the  distribution (\ref{2.8})
belongs to $\EP_{-1}^{-0}$.
The formulas  (\ref{recur}) give the recursion relation:
\begin{equation*}\begin{split}
:\!x(t_1)\dotsm x(t_{n+1})\!:&= \big(:\! x(t_1)\dotsm x(t_n)\!:\,:\!
x(t_{n+1}\!:)\big)
\hat{}-n\big(:\! x(t_1)\dotsm x(t_n)\!: \delta(t_n-t_{n+1})\big)
\hat{}\\
&\quad-n\big(:\!x(t_1)\dotsm x(t_{n-1})\!:1(t_n)\delta(t_n-t_{n+1})
\big)\hat{},\end{split}\end{equation*}
where $\delta=\delta_0$ is the delta function at 0, the index $\hat{} $
stands for the symmetrization of a function.

The second example is the Poisson white noise exponential function
\be{2.9}:\!e^{<x,y>}:\equiv \sum_{n=0}^\infty (n!)^{-1}<\wpx{n},y^{\ot n}>,
\qquad y\in\e'_\C.\end{equation}
If $y\in E_{-p,\C}$, $p>0$, then $:\!e^{<x,y>}\!:\in\EP
_{-p}^{-0}\,$; if $\xi\in E_{p,\C}$, $p\ge0$, then $:\!e^{<x,\xi>}\!:\in
\EP
_{p}^{\kappa}$ with $\kappa<1$, and
\be{2.10}:\!e^{<x,\xi>}\!:\in\EP_{p}^{1}\ \ \ \mbox{if $|\xi|_p<1$}
\end{equation}(we keep the notation $|\cdot|_p$ for the complexified spaces
$E_{p,\C}^{\hot n}$).
We  will return to this function in the next section.

\section{Continuous version theorem and Poisson white noise delta function}
In this section, we will show that every element of the space  $\EP^1$ has a
version
(in the $\LP$-sense) whose restriction to every $E_{-p}$ is continuous on
$E_{-p}$. Such a continuity will be called a continuity on $\e'$ (though
it does not imply the continuity on $\e'$ endowed with the strong dual
topology). This will allow us to introduce a Poisson white noise delta function.
Also, by using the (proof of the) theorem on continuity, we will obtain a theorem on the
explicit form of the Poisson white noise exponential function.

\begin{theorem} Each $\phi$ from the space $\EP^1$ has a version $\tilde\phi$
that is continuous on $\e'$ and is given by
\be{3.1}\tilde\phi(x)=\sum_{n=0}^\infty <\wpx{n},f^{(n)}>,\qquad f^{(n)}\in\e
_\C^{\ho n},\end{equation}
where the Wick powers $\wpx{n}$ are defined by the recursion relation \eqref{recur}
and $\phi$ is the image under the unitary $I$ of the element
$(f^{(n)})_{n=0}^\infty$ of the space $\Gamma_1(\e)$\rom. The series
\eqref{3.1}
converges absolutely and uniformly on every bounded set from $\e'$\rom.
Moreover\rom, $\tilde\phi(x)$ can be extended to the complexification $\e'_\C$
of $\e'$ that $\tilde\phi(z)$ becomes analytic in $\e_\C'$. This extension
is given by the formulas \eqref{3.1} and \eqref{recur} in which $x\in\e'$ is
replaced by
$z\in\e'_\C$\rom. \end{theorem}

\noindent {\it Remark\/} 3.1. The definition of a function analytic in
$\e_\C'$
can be found, e.g., in \cite{Lee,HKPS}.\vs

\noindent{\em Proof}.
Our proof is close in spirit to (part of) the proof of the continuous version
theorem in Gaussian  analysis
proposed by Obata \cite{Ob1,Ob2}. We wish to estimate the norms of $\wpx{n}$
in $E_{-p}^{\hot n}$.

Since $<\tau,\xi^{\ot2}>=<\int_T\delta_t^{\ot2}\,d\nu(t),\xi^{\ot2}>$, we get from
(A.1) and \eqref{2.2} that
\be{3.2}
|\tau|_{-p}\le\int_T|\delta_t|_{-p}^2\,d\nu(t)\le\int_T|\delta_t|^2_{-1}\,d\nu(t)
=\|\delta\|^2,\qquad p\ge1.\end{equation}
By (A.4) and \eqref{starr}, we easily conclude that
\be{3.3}|\DDD^{(n+1)}f^{(n+1)}|_p\le nC_p|f^{(n+1)}|_p,\qquad p\ge1,\end{equation}
i.e., $\DDD^{(n+1)}\colon E_p^{\hot(n+1)}\to E_p^{\hot n}$ is a linear continuous operator with norm
$\le nC_p$.

By \eqref{recur},
\be{3.5}\wpx{n+1}=\wpx{n}\hot\wpx{1}-\DDD^{(n+1)\,*}\wpx{n}-n\wpx{(n-1)}\hot
\tau,\end{equation}
where $\DDD^{(n+1)\,*}\colon \e^{\prime\,\hot n}\to\e^{\hot(n+1)}$ is the dual
operator of $\DDD^{(n+1)}$, and for each $p\ge1$,\linebreak
$\DDD^{(n+1)}\colon E_{-p}^{\hot n}\to E_{-p}^{\hot(n+1)}$ is a continuous operator
with norm $\le nC_p$.

By (A.2) and \eqref{2.2}, we have that $1\in E_{-1}$ and
\be{3.6}|1|_{-p}\le|1|_{-1}\le\int_T|\delta_t|_{-1}\,d\nu(t)\le\|\delta\|_
{\infty},\qquad p\ge1.\end{equation}

Thus, by \eqref{3.2}--\eqref{3.6},
\be{3.7}\begin{split}
|\wpx{(n+1)}|_{-p}&\le|\wpx{n}|_{-p}(|x|_{-p}+\|\delta\|_\infty)\\
&\quad+ nC_p|\wpx{n}|_{-p}+\|\delta\|^2|\wpx{(n+1)}|_{-p},\qquad p\ge1.\end{split}
\end{equation}

Hence, given any fixed $p\ge1$ and $R>0$, we have, for each $x\in\e'$ such that
$|x|_{-p}\le R$,
\begin{gather*}|\wpx{(n+1)}|_{-p}\le nY_{p,R}\max\{|\wpx{n}|_{-p},|\wpx{(n-1)}|
_{-p}\},\\Y_{p,R}=R+\|\delta\|_\infty+C_p+\|\delta\|^2,\end{gather*}
which easily yields that
\be{3.8}
|\wpx{(n+1)}|_{-p}\le n!\,Z^n_{p,R} ,\qquad Z_{p,R}=\max\{1,Y_{p,R}\}.\end{equation}

The assumption (A.3) implies
\be{3.9p}\rho^{mn}|f^{(n)}|_{p+m}\ge|f^{(n)}|_p,\qquad f^{(n)}\in
E_p^{\ho n},\ m\in\n,\end{equation}
and so
\be{3.9} |F^{(n)}|_{-(p+m)}\le\rho^{mn}|F^{(n)}|_{-p},\qquad
F^{(n)}\in E_{-p}^{\ho n},\ m\in\n.\end{equation}

Summing (\ref{3.8}) and (\ref{3.9}) up, we conclude that there exists
$p_1=p_1(p,R)\ge p$ such that $$|\wpx{n}|_{-p_1}\le n!\,2^{-n}.$$
Therefore, for $\phi(x)=\sum_{n=0}^\infty<\wpx{n},f^{(n)}>\in\EP^1$,
we have, for $|x|_{-p}\le R$,
\begin{equation*}\begin{split}
\sum_{n=0}^\infty|<\wpx{n},f^{(n)}>|&\le
\sum_{n=0}^\infty|\wpx{n}|_{-p_1}|f^{(n)}|_{p_1}
\le\sum_{n=0}^\infty n!\,2^{-n}|f^{(n)}|_{p_1}\\
&\le\bigg(\sum_{n=0}^\infty
4^{-n}\bigg)^{1/2}\|\phi\|_{-1,-p_1},\end{split}\end{equation*}
i.e., the series $\sum_{n=0}^\infty<\wpx{n},f^{(n)}>$
converges absolutely and uniformly on every bounded set  in $E_{-p}$
(we recall that every bounded set in $\e'$ endowed with the strong dual topology
is bounded in some space $E_{-p}$ \cite{GeVi}). For any $f^{(n)}\in
\e_{\C}^{\ho n}$,  the function $<\wpx{n},f^{(n)}>$ is a continuous
polynomial of variable $x\in\e'$ (see Proposition~2.2), and so
the function $\tilde\phi(x)$ is continuous on every $E_{-p}$.

Next, we note that if one replaces $x\in\e'$ with $z\in\e'_{\C}$ in
(\ref{recur}) and (\ref{3.1}), then all the above formulas hold true for the
complexified spaces $E_{-p,\C}$.
Since  $<:\!z^{\ot n}\!:,f^{(n)}>$, $f^{(n)}\in\e_{\C}^{\ho n}$, is a
continuous polynomial of variable $z$, it is an analytic function in
every $E_{-p,\C}$. The series
$$\tilde\phi(z)=\sum_{n=0}^\infty<:\!z^{\ot n}\!:,f^{(n)}>$$
converges absolutely and uniformly on every bounded set in $E_{-p,\C}$,
so $\tilde\phi(z)$ is analytic in every $E_{-p,\C}$. By \cite{Lee},
Theorem~A.2, we conclude that $\tilde\phi(z)$ is analytic in $\e'_{\C}$.
\hfill $\Box$\vs

Having obtained  the continuous version theorem, we are able now
to define a Poisson white noise delta function. So,  we put, for
each $y\in\e'$,
\[ \tilde\delta_y=\sum_{n=0}^\infty
<\x{n},(n!)^{-1}\y{n}>.\]

\begin{prop} For each $y\in\e'$, $\tilde\delta_y$ belongs to $\EP^{-1}$\rom,
and for each $\phi\in\EP^1$\rom,
\be{3.10} \ll\tilde\delta_y,\phi\gg=\tilde\phi(y),\end{equation}
where $\tilde\phi$ is the continuous version of $\phi$ defined in
Theorem~\rom{3.1.}
Moreover\rom, the following mapping is continuous\rom:
$$\e'\ni y\to\tilde\delta_y\in\EP^{-1}.$$\end{prop}

\noindent{\em Proof}. Let us fix $y\in\e'$, and let $p>0$ be such
that $y\in E_{-p}$,
$|y|_{-p}=R$. Let $p_1\ge p$ be chosen so that
\be{3.11} (n!)^{-1}|\y n|_{-p_1}\le 2^{-n}\end{equation}
(see the proof of Theorem~3.1). Then
$$\|\tilde\delta_y\|_{-1,-p_1}^2=\sum_{n=0}^\infty (n!)^{-2}|\y n|
_{-p_1}^2\le\sum_{n=0}^\infty 4^{-n}<\infty,$$
whence $\tilde\delta_y\in\EP^{-1}_{-p_1}$ and (\ref{3.10}) holds.

Let a sequence $\{y_j\,|\, j\in\n\}\in E_{-p}$ tend to a $y$ in $E_{-p}$.
Then, there is $R>0$ such that $|y_j|_{-p}\le R$ for all $j\in\n$, and so
$|y|_{-p}\le R$. Choose $p_1$ so that  (\ref{3.11}) holds for any $y\in
E_{-p}$ with $|y|_{-p}\le R$. Let us show that
\be{3.12} \tilde\delta_{y_j}\to \tilde\delta_y\ \ \mbox{in
$\EP_{-p_1}^{-1}$ as $j\to\infty$}.\end{equation}
As follows from the proof  of Theorem~3.1,
$$\yj{n}\to\y{n}\ \ \mbox{in $E_{-p_1}^{\hot n}$, $n\in\Z_+$},$$
whence
\begin{equation}\label{3.13}
 <\x{n},(n!)^{-1}\yj{n}>\to
<\x{n},(n!)^{-1}\y{n}>\ \ \mbox{in $\EP_{-p_1}^{-1}$ as $j\to\infty$,
$n\in\Z_+$.}\end{equation}
Next, for any $n,m\in\n$, $n>m$, we get by (\ref{3.11})
\be{3.14}\begin{split}
\sum_{i=m}^n\|<\x{i},(i!)^{-1}(\yj{i}-\y{i})>\|^2_{-1,-p_1}&\le
\sum_{i=m}^n2(i!)^{-2}\big(|\yj{i}|_{-p_1}^2
+|\y{i}|_{-p_1}^4\big)\\&\le
4\sum_{i=m}^n4^{-i}.\end{split}\end{equation}
From \eqref{3.13} and \eqref{3.14},  we easily conclude that (\ref{3.12})
holds.\hfill $\Box$\vs

Now, we will prove a theorem on the evident form
of the Poisson white noise exponential function. This theorem is a refinement
of a corresponding result of \cite{ItKu}.

\begin{theorem} We have
\be{3.16}\begin{gathered}
:\!e^{<x,\xi>}\!:=\exp\Big[ <x,\log(1+\xi)>-\int_T\xi(t)\,d\nu(t)\Big],\\
\xi\in\e_\C,\ |\xi|_1<\max\{1,C_1\},\end{gathered}\end{equation}
which holds for $\muP$-a\rom.a\rom.\ $x\in\e'$\rom, more exactly\rom,
for all $x\in E_{-1}$\rom,
which is a set of full $\muP$ measure\rom.\end{theorem}

\noindent{\em Proof}. We divide the proof into steps.

1. Let us fix an arbitrary set $\alpha\subset T$ of finite $\nu$ measure and
put
$$f(t)=\sum_{j=1}^k\omega_j\chi_{\alpha_j}(t),\qquad \omega_j\in\C,$$
where $\alpha_j$, $j=1,\dots,k$, are disjoint subsets of $\alpha$,
$\bigcup_{j=1}^k\alpha_j=\alpha$. Then, by Proposition~2.1,
\bea
\ds <\x{n},\chi_{\alpha_1}^{\ot n_1}\hot\dotsm\hot\chi_{\alpha_k}
^{\ot n_k}>=\prod_{j=1}^k<\x{n_j},\chi_{\alpha_j}^{\ot n_j}>,\quad
n_1+\dotsb+n_k=n,&\nonumber\\
\ds<\x{n_j},\chi_{\alpha_j}^{\ot n_j}>=C_{n_j}(<x,\chi_{\alpha_j}>,
\nu(\alpha_j)),\label{dyr}&\end{eqnarray}
where $C_n(u,\tau)$ are the Charlier polynomials with parameter $\tau$:
\be{dyrr}\sum_{n=0}^\infty\frac{\omega^n}{n!}\,C_n(u,\tau)=\exp\big[u\log(1+\omega)
-\omega\tau\big].\end{equation}
 Therefore, by \eqref{dyr}, \eqref{dyrr}, we have (cf.\ \cite{It,ItKu})
\be{3.15}\begin{split}
&\sum_{n=0}^\infty (n!)^{-1}<\x{n}, f^{\ot n}>\\
&\qquad =\exp\bigg[<x,\sum_{j=1}^k\chi_{\alpha_j}\log(1+\omega_j)>-\int_T
\Big(\sum_{j=1}^k \omega_j\chi_{\alpha_j}(t)\Big)d\nu(t)\bigg]\\
&\qquad=\exp\bigg[<x,\chi_\alpha\log(1+f)>-\int_Tf(t)\,d\nu(t)\bigg],
\end{split}\end{equation}
the equalities making sense for $|\omega_j|<1$.

2. For any $n\in\n$, let us consider the following triple
$$\big(\hat L^\infty(\alpha^n,\nu^n)\big)'\supset\hat L^2(\alpha^n,\nu^n)
\supset \hat L^\infty(\alpha^n,\nu^n),$$
where $\hat L^\infty(\alpha^n,\nu^n)$ is the subspace of
$L^\infty(\alpha^n,\nu^n)$ consisting of symmetric functions on $\alpha^n$,
and $ \big(\hat L^\infty(\alpha^n,\nu^n)\big)'$
is the dual of $\hat L^\infty(T^n,\nu^n)$ with respect to zero space
$\hat L^2(\alpha^n,\nu^n)$. Let us fix an arbitrary $x=x(t)\in
L^1(\alpha,\nu)$.  We suppose that $\x{n}$ are defined
by the recursion relation (\ref{recur}) and these $\x{n}$ are understood as
elements of $\big(\hat L^\infty(\alpha^n,\nu^n)
\big)'$. By analogy with the proof of Theorem~3.1, the following
estimate can be proved
$$\|\x{n}\|_{\big(\hat L^\infty(\alpha^n,\nu^n)\big)'}\le n!\,V^n,$$
where $V$ is a positive constant, depending on $\alpha$ and $x$. Therefore,
for any $f\in L^\infty(\alpha,\nu)$, we have
$$\sum_{n=0}^\infty (n!)^{-1}|<\x{n},f^{\ot n}>|\le\sum_{n=0}^\infty
\big(V\|f\|_{L^\infty(\alpha,\nu)}\big)^n.$$
Hence, $x\in L^1(\alpha,\nu)$
being fixed,
$$:\!e^{<x,f>}\!:=\sum_{n=0}^\infty (n!)^{-1}<\x{n},f^{\ot n}>$$
is a continuous function of variable $f\in L^\infty(\alpha,\nu)$
in the ball in $L^\infty(\alpha,\nu)$ of radius $(2V)^{-1}$ centered at
0.

Let us fix $\xi\in\e_\C$ such that
$|\xi|_1\le(2V\|\delta\|_\infty)^{-1}$,
and so, by (A.2),
$\|\xi\|_{L^\infty(T,\nu)}\le(2V)^{-1}$.
Approximate the function $\xi\chi_\alpha$ in the $L^\infty(T,\nu)$
norm by step functions. Since (\ref{3.15}) holds for all $x\in L^1(T,\nu)$,
we conclude that (\ref{3.15}) holds true for our fixed $x$ if $f$ is
replaced by $\xi\chi_\alpha$.

3. Let us denote by $f(x;\xi)$ the function on the right hand side of
\eqref{3.16}. By virtue of (A.4), $\log(1+\xi)=\sum_{n=1}^\infty\frac{(-1)^
{n+1}}{n}\,\xi^{\ot n}\in E_{1,\C}$
provided $\xi\in\e$ and $|\xi|_1<R=\max\{1,C_1\}$.
Thus, taking also to notice (A.2), we conclude that $f(x;\xi)$ is well defined
for all $x\in E_{-1}$ and the above $\xi$. Moreover, for each $x\in E_{-1}$
fixed, $f(x;\xi)$ as a function of $\xi\in E_{1,\C}$ is analytic in the ball $
|\xi|_1<R$. Therefore, see, e.g., the proof of
Theorem~3 in \cite{KLS}, it can be represented in the form
$$f(x;\xi)=\sum_{n=0}^\infty (n!)^{-1}<F^{(n)}(x),\xi^{\ot n}>, \qquad
F^{(n)}(x)\in\e_{\C}^{\prime\, \ho n}.$$
Moreover, since $f(x;\xi)\in\R$ for any $\xi\in\e$,
we conclude that actually $F^{(n)}(x)\in\e^{\prime\,\ho n}$, i.e.,
$<F^{(n)}(x),f^{(n)}>\in\R$ for any $f^{(n)}\in\e^{\ho n}$.

As follows from 1--2, for any $x\in L^1(T,\nu)$ and any set $\alpha$
of finite $\nu$ measure, we have
$$:\!e^{<x,\xi>}\!:=\sum_{n=0}^\infty (n!)^{-1}<\x{n},\xi^{\ot n}>=
\sum_{n=0}^\infty(n!)^{-1}<F^{(n)}(x),\xi^{\ot n}>=f(x;\xi),$$
where $\xi\in\e$ is such that it vanishes outside of $\alpha$ and $|\xi|
_1\le{\mathcal Y}_\alpha$, ${\mathcal Y}_\alpha$ a positive constant determined by
$x$ and $\alpha$.  From here, because of (A.5), $\x n=F^{(n)}(x)$. Thus, (\ref{3.16})
holds for all $x\in L^1(T,\nu)$.

4. It  follows from the proof of Theorem~3.1 that there exist $p_0\ge1$
and $r>0$ such that, for any fixed $\xi\in\e$, $|\xi|_{p_0}\le r$,
$:\!e^{<x,\xi>}\!:$ is a continuous function of variable $x\in E_{-1}$.
On the other hand, $f(x;\xi)$ is continuous on $E_{-1}$ for all
$\xi\in\e$, $|\xi|_1<R$. Let us fix an arbitrary $\xi\in\e$,
$|\xi|_{p_0}\le r$. We know that
\be{3.17}
:\!e^{<x,\xi>}\!:=f(x;\xi)\quad \mbox{for all $x\in L^1(T,\nu)$}\end{equation}
(we assume that $r\le R$). Extending (\ref{3.17}) by continuity
(since $\e\subset L^1(T,\nu)$ and $\e$ is a dense subset of $E_{-1}$,
$L^1(T,\nu)$ is also a dense subset of $E_{-1}$), we conclude that
(\ref{3.17}) holds true for all $x\in E_{-1}$ and $\xi\in\e$, $|\xi|_{p_0}\le
r$. At last, analogously to 3, we get the conclusion of the theorem.
\hfill $\Box$\vs

\noindent{\it Remark\/} 3.2. It follows from the proof of Theorem~3.2
that the formula (\ref{3.16}) holds for all $x\in E_{-p}$, $p\ge1$,
provided $\xi\in \e_{\C}$ and $|\xi|_p<\big(\max\{1,C_p\}\big)^
{-1}$.

\begin{corollary} For each $y\in\e'$\rom, $y\ne0$\rom, the Poisson white noise
delta function $\tilde \delta_y$ does not belong  to any space $\EP^\kappa$
with $\kappa<1$\rom.\end{corollary}

\noindent{\it Remark\/} 3.3. This statement shows that there is no sense in
trying
to prove the continuous version theorem  for any $\EP^{-\kappa}$ with
$\kappa<1$.\vs

\noindent{\it Proof}. Following \cite{ItKu,KS}, for any $\Phi\in\EP^{-\kappa}$
with $\kappa<1$, we define the $\mathcal S_{\text{P}}$-transform of $\Phi$ by
\[\mathcal S_{\text{P}}[\Phi](\xi)=\ll\Phi,:\!e^{<x,\xi>}\!:\gg,\qquad\xi\in\e_\C.\]
Since the set of $:\!e^{<x,\xi>}\!:$ is total in each $\EP^{-\kappa}$, the
$\mathcal S_{\text{P}}$-transform   uniquely defines $\Phi$. Moreover, for each $
\Phi\in\EP^{-1}$, there is $p>0$ such that $\Phi\in\EP_{-p}^{-1}$, and we set,
taking to notice \eqref{2.10},
\be{3.18}\mathcal S_{\text{P}}[\Phi](\xi)
=\ll\Phi,:\!e^{<x,\xi>}\!:\gg,\qquad\xi\in\e_\C,\ |\xi|_p<1.\end{equation}
As follows from  \cite{KLS}, the $\mathcal S_{\text{P}}$-transform defined by
\eqref{3.18}  uniquely determines $\Phi\in\EP^{-1}$.

Since all the above constructions are totally isomorphic to the Gaussian case,
all the characterization  theorems of the spaces $\EP^{\kappa}$,  $\kappa\in
[-1,1]$, in terms of their $\mathcal S_{\text{p}}$-transforms hold true, e.g.,
\cite{BeKo,HKPS,Ob2,KS,KLS}.

Let us fix $y\in\e'$, $y\ne0$, and let $p>0$ be such that $y\in E_{-p}$.
By virtue of Proposition 3.1, we know that {\it a priori\/} $\tilde
\delta_y\in\EP^{-1}$. By Theorem~3.2 (more exactly, by Remark~3.2), we get
that $\mathcal S_{\text{P}}[\tilde \delta_y](\xi)$ is equal to the right
hand side of \eqref{3.16} provided $|\xi|_p<\max\{1,C_p\}^{-1}$.
But the function $\mathcal S_{\text{P}}[\tilde \delta_y](\xi)$ is analytic
only in a neighborhood of zero in $\e_\C$, but not in the whole $\e_\C$.
Therefore, by the characterization theorems, which state particularly
that the $\mathcal S_{\text P}$-transform of an element of $\EP^{-\kappa}$ with
$\kappa<1$ is a function analytic in $\e_\C$, we obtain the desired statement.
\hfill$\Box$

\begin{corollary} For each $t\in T$\rom, define a linear continuous operator
$\di_t$ in $\EP^{1}$ by
\[\di_t<\x{n},f^{(n)}>=n<\x{(n-1)},f^{(n)}(\cdot,\dots,\cdot,t)>,\]
$\di_t$ is called the operator of Hida differentiation at $t$\rom. Then\rom, for each
$\phi\in\EP^{1},$
\[\big(\di_t\phi\big)^\sim(x)=\tilde\phi(x+\delta_t)-\tilde\phi(x),\qquad
x\in\e',\]
where $\tilde\phi$ denotes the continuous version of $\phi$ defined in Theorem~\rom{3.1.}
\end{corollary}

\noindent{\it Proof}. For $\phi\in\mathcal P(\e')$, the corollary was proved in
\cite{It,ItKu}. For an arbitrary $\phi\in\EP^1$, let us choose a sequence
$(\phi_j)_{j=1}^\infty\in\mathcal P(\e')$ such that $\phi_j\to\phi$ in $\EP^1$.
Then, by Proposition~3.1,
\be{3.19}
\tilde\phi_j(x+\delta_t)-\tilde \phi_j(x)=\ll\tilde \delta_{x+\delta_t}
-\tilde\delta_x,\phi_j\gg
\to\ll\tilde\delta_{x+\delta_t}-\tilde \delta_x,\phi\gg=\tilde\phi(x+\delta_t)
-\tilde\phi(x).
\end{equation}
On the other hand,
\be{3.20}
\tilde\phi_j(x+\delta_t)-\tilde\phi_j(x)=(\di_t\phi_j)^\sim(x)=\ll\tilde\delta_x,
\di_t\phi_j\gg
\to\ll\tilde \delta_x,\di_t\phi\gg=(\di_t\phi)^\sim(x).
\end{equation}
Combining \eqref{3.19} and \eqref{3.20} gives the corolary.\hfill$\Box$\vs

\noindent{\it Remark\/} 3.4. It is worth to compare Corollary 3.2 with the
results of Nualart and Vives \cite{NuVi}

\section{Multiplication in $\EP^1$}
\begin{theorem}
The space $\EP^1$ is an algebra under pointwise multiplication of functions\rom.
More exactly\rom, for any $\phi,\psi\in\EP^1$ and $p\ge1$\rom, there is
$\operatorname{const}>0$ such that
\[
\|\phi\psi\|_{1,p}\le\operatorname{const}\|\phi\|_{1,p+1}\|\psi\|_{1.p+q},\]
where $q\in\n$ is chosen so that
\[\rho^q<(1-\rho)^2Y_p^{-1},\qquad Y_p=\max\{1,C_p\}\max\{1,\|\delta\|
_\infty\},\]
$\rho$ the constant from \rom{(A.3).} In particular\rom, for any $\phi
\in\EP^1$\rom, the operator of multiplication by $\phi$
acts continuously from $\EP^1_{p+1}$ into
$\EP_p^1$ for each $p\ge1.$\end{theorem}

\noindent{\it Proof}. The proof is rather analogous to that of
Proposition~6.5 in \cite{ItKu}, so we  only note some new points.

Let $f^{(n)}\in\e^{\hot n}$, $n\in\Z_+$. Then, for the operator $A_{j^*,k,j}
(f^{(n)})$, $j^*+k+j=n$, defined in \cite{ItKu}, we have the estimate (compare
with Proposition~4.5 in \cite{ItKu}):
\be{4.1}\|A_{j^*,k,j}(f^{(n)})\phi\|_{1,p}\le\rho^{(j^*+k)}(j^*+k)!\,
(1-\rho)^{-(j^*+k+1)}C_p^{k+j}\|\delta\|_\infty^j\rho^{j(p-1)}|f^{(n)}|_p
\|\phi\|_{1,p+1},\end{equation}
so that $A_{j^*,k,j}(f^{(n)})$ is a continuous operator in $\EP^1$.

Let $\psi\in\EP^1$ be of the form  $\psi(x)=\sum_{n=0}^\infty<\x{n},f^{(n)}>$.
Then \cite{ItKu}
\be{4.2}<\x{n},f^{(n)}>\phi(x)=\sum_{j^*+k+j=n}\frac{n!}{j^*!\,k!\,j!}\,
A_{j^*,k,j}(f^{(n)})\phi.\end{equation}
By \eqref{4.1} and \eqref{4.2}
\be{4.3}\|\phi\psi\|_{1,p}\le\|\phi\|_{1,p+1}\sum_{n=0}^\infty n!\, (1-\rho)
^{-(n-1)}Y_p^n|f^{(n)}|_p\sum_{j^*+k+j=n}\frac{(j^*+k)!}{j^*!\,k!\,j!}\,
\rho^k.\end{equation}
By using the estimate (3.25) in \cite{ItKu}, we have
\[\frac{1}{j^*!\,j!}\,\sup_{k\ge0}\frac{(j^*+k)!}{k!}\,\rho^k\le\frac{1}{j!}\,
(1-\rho)^{-j^*-1},\]
which, upon \eqref{4.3} and \eqref{3.9p}, gives
\begin{equation*}\begin{split}
\|\phi\psi\|_{1,p}&\le\|\phi\|_{1,p+1}\sum_{n=0}^\infty n!\,(1-\rho)^{-2(n+1)}
Y_p^n\rho^{nq}|f^{(n)}|_{p+q}\,\tfrac12(n+1)(n+2)\\
&\le\tfrac12\|\phi\|_{1,p+1}\|\psi\|_{1,p+q}(1-\rho)^{-2}\bigg(\sum_{n=0}^\infty
(n+1)^2(n+2)^2\big[(1-\rho)^{-2}Y_p\rho^q\big]^{2n}\bigg)^{1/2},\end{split}\end{equation*}
which gives the theorem.\hfill$\Box$\vs

\noindent {\it Remark\/} 4.1. By analogy with the proof of Theorem 4.1, one can easily
verify that each space $\EP^{\kappa}$ with $\kappa>1$ has an algebraic
structure.\vs

Following  \cite{ItKu}, for each $t\in T$, we define the operator of
Poisson coordinate multiplication
\[x(t){\cdot}=(\di_t^*+1)(\di_t-1)=\di_t^*\di_t+\di_t+\di_t^*+1,\]
where $\di_t^*\colon \EP^{-1}\to\EP^{-1}$ is the dual of the operator
$\di_t\colon\EP^1\to\EP^1$ defined in Corollary~3.2. Evidently, $x(t){\cdot}$
is a continuous operator from $\EP^1$ to $\EP^{-1}$.

\begin{corollary}For $\phi,\psi\in\EP^1$ and $t\in T$\rom,
\[\ll x(t)\cdot\phi,\psi\gg=\ll x(t),\overline\phi\psi\gg,\]
where $\overline \phi$ is the complex conjugate of $\phi$\rom. Here\rom,
$x(t)=<\x1,\delta_t>+1\in\EP^{-1}.$\end{corollary}

\noindent{\it Proof}. Let us choose an arbitrary sequence $(\xi_j)_{j=1}
^\infty\subset \e$ such that $\xi_j\to\delta_t$ in $\e'$ (notice that
$\e$ is dense in $\e'$). Denote by $<x,\xi_j>{\cdot}$ the operator of multiplication by
the function $<x,\xi>$. Evidently,
\[\ll<x,\xi_j>\cdot\phi,\psi\gg=\ll<x,\xi>,\overline\phi\psi\gg.\]

It remains only to note that $<x,\xi_j>\to x(t)$ in $\EP^{-1}$ and
$<x,\xi_j>\cdot\phi\to x(t)\cdot\phi$ in $\EP^{-1}$ (see \cite{ItKu,Ly}).
\hfill$\Box$\vs

\begin{center}{\bf Acknowledgments}\end{center}
The author is grateful to Prof.\ Yu.\ M. Berezansky for encouraging him
to study non-Gaussian analysis. The author also expresses his graduate
to Profs.~Yu.~G.~Kondratiev  and G.~F.~Us for useful discussions.

\noindent Institute of Mathematics of Ukrainian National Academy of Sciences\vs

\noindent Institute of Mathematics\\
Ukrainian National Academy of Sciences\\
3 Tereshchenkivska St., Kiev 252601\\
Ukraine
\end{document}